\documentclass{amsart}%
\usepackage[english]{babel}
\usepackage{amsthm}
\usepackage[letterpaper,top=2cm,bottom=2cm,left=3cm,right=3cm,marginparwidth=1.75cm]{geometry}

\usepackage{amsmath}
\usepackage{mathtools}
\usepackage{graphicx}
\usepackage[colorlinks=true, allcolors=blue]{hyperref}
\usepackage{algorithm,algorithmicx,algpseudocode}%
\usepackage{amssymb}
\title{Anisotropy of quadratic forms over global fields of characteristic $\neq$ 2 is Diophantine}
\author{Guang Hu}
\address{The Mathematics School, Shandong University, 27 Shanda Nanlu, Jinan, Shandong, P. R. China}
\address{Faculty of Mathematics, Technische Universität Dresden, 01062 Dresden, Germany}
\email{202490000031@sdu.edu.cn,huguang@mail.ustc.edu.cn,guang.hu@mailbox.tu-dresden.de}
\begin{document}

\def\Br{\textup{Br}}
\def\H{\textup{H}}
\def\Sp{\textup{Spec}}
\def\inv{\textup{inv}}
\def\BQ{\mathbb{Q}}
\def\BR{\mathbb{R}}
\def\BZ{\mathbb{Z}}
\def\BN{\mathbb{N}}
\def\BA{\mathbb{A}}
\def\FA{\mathbf{A}}
\def\BP{\mathbb{P}}
\def\BG{\mathbb{G}}
\def\BF{\mathbb{F}}
\def\BC{\mathbb{C}}
\def\CO{\mathcal{O}}
\def\KB{\mathfrak{B}}
\def\KT{\mathfrak{T}}
\def\Km{\mathfrak{m}}
\def\KS{\mathfrak{S}}
\def\Kp{\mathfrak{p}}
\def\Kq{\mathfrak{q}}
\def\KR{\mathfrak{R}}
\def\N{\textup{N}}
\def\NS{\textup{NS}}
\def\char{\textup{char}}
\def\rk{\textup{rk}}
\def\Pic{\textup{Pic}}
\def\Hom{\textup{Hom}}
\def\GL{\textup{GL}}
\def\Gal{\textup{Gal}}
\def\PGL{\textup{PGL}}
\def\,{\textup{，}}
\def\Gal{\textup{Gal}}
\def\M{\textup{M}}
\def\Az{\textup{Az}}
\def\.{\textup{。}}
\def\inj{\hookrightarrow}
\def\ev{\textup{ev}}
\def\et{\textup{ét}}
\newtheorem{thm}{Theorem}[section]
\newtheorem{cor}[thm]{Corollary}
\newtheorem{lem}[thm]{Lemma}
\newtheorem{prob}[thm]{Problem}
\newtheorem{pty}[thm]{Proposition}
\theoremstyle{definition}
\newtheorem{defin}[thm]{Definition}
\newtheorem{rem}[thm]{Remark}
\newtheorem{exa}[thm]{Example}
\begin{abstract}
We prove that the set of anisotropic quadratic forms over global fields of characteristic different from 2 is a diophantine set. Our proof builds upon and extends the method of Koenigsmann, using tools from class field theory, the local-global principle, and advances on the diophantine definability of non-norm sets over global fields.
\end{abstract}
\maketitle

\section{Introduction}
Let $K$ be a global field. Roughly speaking, a diophantine set is a set that can be described by the solvability of a family of polynomial equations. Equivalently, we give the following definition.
\begin{defin}\label{d1}
A subset $S \subset K^n$ is called a diophantine set over $K$ if there exists a polynomial
\[
f(a_1,\ldots,a_n; x_1,\ldots,x_m) \in K[a_1,\ldots,a_n; x_1,\ldots,x_m]
\]
such that for each $(a_1,\ldots,a_n) \in K^n$ we have
\[
(a_1,\ldots,a_n) \in S \iff \exists (x_1,\ldots,x_m)\in K^m \text{ such that } f(a_1,\ldots,a_n; x_1,\ldots,x_m) = 0.
\]
\end{defin}
The study of diophantine sets is closely connected with Hilbert’s Tenth Problem. Matiyasevich proved that the problem admits a negative answer over $\mathbb{Z}$ (see \cite{mt}). Furthermore, if $\mathbb{Z}$ is diophantine over $\mathbb{Q}$, then Hilbert’s Tenth Problem over $\mathbb{Q}$ would also have a negative answer. However, this remains an open problem.

Besides their connection to foundational questions in logic and number theory, the study of diophantine sets is of intrinsic interest as an inverse to the problem of solving polynomial equations. While the classical problem asks for the solutions of given polynomial equations, the study of diophantine sets addresses the reverse question: which subsets of $K^n$ can be realized as the set of solutions to some family of polynomial equations.

A significant advance in this direction was achieved by Poonen \cite{Poo09}, who proved, using the Brauer--Manin obstruction, that the set of non-squares in number fields is diophantine. This provided one of the first nontrivial examples of such sets.

This result has been generalized to various situations in recent years. Colliot-Thélène and Jan Van Geel \cite{CTGV} extended Poonen’s result to the complements of $n$-th powers over global fields. Koenigsmann \cite{10.2307/24735167} introduced powerful new methods for constructing diophantine sets over $\mathbb{Q}$, notably the set of non-norms for quadratic extensions.

More generally, for a quadratic extension $L/K$ of a field $K$, the set of non-norms 
\[
\{x \in K^{\times} \mid x \notin N_{L/K}(L^{\times})\}
\]
is diophantine. In fact, a stronger statement holds: 
\[
\{(x,y)\in (K^{\times})^2 \mid x \notin N_{K(\sqrt{y})/K}(K(\sqrt{y})^{\times})\}
\]
is diophantine. Park \cite{JP} extended Koenigsmann’s method to number fields, and Eisenträger and Morrison \cite{KTU} later applied it to global function fields of odd characteristic. Morrison \cite{Morrison2017diophantineDO} further generalized the result to cyclic extensions of degree prime to the characteristic.

Another result was obtained by Philip Dittmann \cite{DIT}. He generalized the method uniformly to apply to any global field and proved that the set of coefficients of univariate polynomials that have no solutions in a global field is diophantine. As a corollary, he showed that the set of coefficients of irreducible multivariate polynomials is diophantine.

In this paper, we study diophantine sets defined by quadratic forms. We begin by recalling the notions of isotropy and anisotropy of quadratic forms.
\begin{defin}
Let $F$ be a field and let $f$ be a quadratic form in $m$ variables over $F$.
We say that $f$ is \emph{isotropic over $F$} if there exists a tuple
$(x_1,\dots,x_m) \in F^m$, not all zero, such that
\[
f(x_1,\dots,x_m)=0.
\]
Otherwise, we say that $f$ is \emph{anisotropic over $F$}.
\end{defin}

Our goal is to show that the set of coefficients of anisotropic quadratic forms over a global field of characteristic not equal to $2$ is a diophantine set. In any case, the set of coefficients of isotropic quadratic forms is obviously diophantine. So, the set we deal with here is the complement of this obvious diophantine set, just like the works we mentioned above. Our main theorem is the following.
\begin{thm}\label{t1}
Let $K$ be a global field with $\textup{char}(K)\neq 2$ and let $m\in\BN$. Then
\[
\bigl\{(a_i)\in (K^{\times})^{m} : \sum_{i=1}^m a_i x_i^2 \text{ is anisotropic over } K \bigr\}
\]
forms a diophantine set over $K$.
\end{thm}
\begin{rem}
In the notation $(K^{\times})^{m}$ we mean the Cartesian product of $m$ copies of $K^{\times}$, 
i.e.\ the set of $m$-tuples of elements of $K^{\times}$, and not the set of $m$th powers in $K^{\times}$, 
which we denote by ${K^{\times}}^{m}$.
\end{rem}

As a corollary of Theorem \ref{t1}, we obtain the following.

\begin{cor}\label{t2}
Let $K$ be a global field with $\mathrm{char}(K)\neq 2$, and let $m\in \BN$. Then the set
\[
\bigl\{(a_{ij})\in (K)^{m^2} \;\big|\;
\sum_{1\leq i,j\leq m} a_{ij} x_i x_j
\text{ is anisotropic over } K \bigr\}
\]
is a diophantine set over $K$.
\end{cor}
To prove Theorem~\ref{t1}, we distinguish cases according to the value of $m$.
For $m \geq 5$, the result follows from local–global principles.
For $m = 1, 2,$ or $3$, the result is already known, in the context of ``non-squares'' and ``non-norms''.
The only genuinely new case is $m = 4$, for which we employ generalizations of Koenigsmann's method developed in \cite{DIT, JP, KTU}.
\section{Preliminaries}

From now on, we denote by $S_K$ the set of primes of $K$, including both finite and infinite primes.
By a prime of $K$ we mean an equivalence class of nontrivial absolute values on $K$.

If $K$ is a number field, we write $S_K^\infty$ for the set of infinite primes of $K$.
For $\mathfrak{p} \in S_K^\infty$, if the corresponding completion $K_{\mathfrak{p}}$ is isomorphic to the real field $\mathbb{R}$, we call $\mathfrak{p}$ a \emph{real infinite prime};
if $K_{\mathfrak{p}} \simeq \mathbb{C}$, then $\mathfrak{p}$ is called a \emph{complex infinite prime}.
If $\mathfrak{p}$ is finite, we also use $\mathfrak{p}$ to denote the corresponding prime ideal, and we write $v_{\mathfrak{p}}$ for the associated normalized discrete valuation on $K$.
Let
\[
S_K^{f} \coloneqq S_K \setminus S_K^{\infty}
\]
denote the set of finite primes of $K$.

For $a,b \in K_{\mathfrak{p}}^\times$, the notation $(a,b)_{\mathfrak{p}}$ denotes the Hilbert symbol at $\mathfrak{p}$.
Note that the Hilbert symbol $(a,b)_{\mathfrak{p}}$ defines a nondegenerate bilinear pairing on
$K_{\mathfrak{p}}^\times / K_{\mathfrak{p}}^{\times 2}$.

For $a \in K_{\mathfrak{p}}^\times / K_{\mathfrak{p}}^{\times 2}$, we define the group homomorphism
\[
T_a \colon K_{\mathfrak{p}}^\times / K_{\mathfrak{p}}^{\times 2} \longrightarrow \{\pm 1\}
\]
by
\[
x \longmapsto (x,a)_{\mathfrak{p}},
\]
where we regard $\{\pm 1\}$ as a multiplicative group.
Let $f$ be a quadratic form in $m$ variables over $K$, and let $a \in K^\times$. We say that $f$ \textit{represents} $a$ over $K$ if there exists a tuple $(x_1, \dots, x_m) \in K^m$ such that 
\[
f(x_1, \dots, x_m) = a.
\]
We say that $f$ is \textit{isotropic} (or \textit{anisotropic}) over $\Kp$ if it is isotropic (or anisotropic) over the completion $K_{\Kp}$.

\begin{lem} \label{half}
With the notations above, let $\epsilon, \epsilon' \in \{\pm 1\}$ and let $a, a' \in K_{\mathfrak{p}}^\times / K_{\mathfrak{p}}^{\times 2}$ be nontrivial. Then
\[
T_a^{-1}(\epsilon) \cap T_{a'}^{-1}(\epsilon') = \emptyset
\quad \Longleftrightarrow \quad
a = a' \text{ and } \epsilon = -\epsilon'.
\]
\end{lem}
\begin{proof}
$\Leftarrow:$ This is obvious.\\

$\Rightarrow:$ By 63:9 of \cite{iqf}, the index $[K_{\mathfrak{p}}^{\times} : K_{\mathfrak{p}}^{\times 2}]$ is finite. 
Since $a \neq 1$ and the Hilbert symbol is nondegenerate, $T_a$ is nontrivial. 
Hence, the cardinalities of $T_a^{-1}(1)$ and $T_a^{-1}(-1)$ are both equal to half of $[K_{\mathfrak{p}}^{\times} : K_{\mathfrak{p}}^{\times 2}]$.

For $\epsilon, \epsilon' \in \{\pm 1\}$ and nontrivial $a, a' \in K_{\mathfrak{p}}^{\times} / K_{\mathfrak{p}}^{\times 2}$, we have 
\[
T_a^{-1}(\epsilon) \cap T_{a'}^{-1}(\epsilon') = \emptyset
\quad \Longleftrightarrow \quad
T_a^{-1}(\epsilon) \sqcup T_{a'}^{-1}(\epsilon') = K_{\mathfrak{p}}^{\times} / K_{\mathfrak{p}}^{\times 2},
\]
by a cardinality argument.

The complement of $T_a^{-1}(1)$ in $K_{\mathfrak{p}}^{\times} / K_{\mathfrak{p}}^{\times 2}$ is $T_a^{-1}(-1)$. 
Hence,
\[
T_a^{-1}(\epsilon) \sqcup T_{a'}^{-1}(\epsilon') = K_{\mathfrak{p}}^{\times} / K_{\mathfrak{p}}^{\times 2}
\quad \text{is equivalent to} \quad
T_a^{-1}(1) = T_{a'}^{-1}(1) \text{ or } T_a^{-1}(1) = T_{a'}^{-1}(-1).
\]

However, $T_a^{-1}(1) = T_{a'}^{-1}(-1)$ is impossible, since this would imply
\[
(x, aa')_{\mathfrak{p}} = (x,a)_{\mathfrak{p}} (x,a')_{\mathfrak{p}} = -1 \quad \forall x \in K_{\mathfrak{p}}^{\times} / K_{\mathfrak{p}}^{\times 2},
\]
which cannot happen.

Also, it is easy to verify that
\begin{align*}
T_a^{-1}(1) = T_{a'}^{-1}(1) 
&\Longleftrightarrow (x,a)_{\mathfrak{p}} = (x,a')_{\mathfrak{p}} \quad \forall x \in K_{\mathfrak{p}}^{\times} / K_{\mathfrak{p}}^{\times 2} \\
&\Longleftrightarrow (x, aa')_{\mathfrak{p}} = 1 \quad \forall x \in K_{\mathfrak{p}}^{\times} / K_{\mathfrak{p}}^{\times 2} \\
&\Longleftrightarrow T_{aa'} \text{ is trivial} \\
&\Longleftrightarrow aa' = 1 \\
&\Longleftrightarrow a = a'.
\end{align*}

In conclusion, we have
\[
T_a^{-1}(\epsilon) \cap T_{a'}^{-1}(\epsilon') = \emptyset \quad \Longleftrightarrow \quad a = a' \text{ and } \epsilon = -\epsilon'.
\]
\end{proof}

\subsection{Quadratic forms}
The arithmetic of quadratic forms over a global field is well understood. 
When the characteristic of the field is not $2$, the situation simplifies further, as the local-global principle applies.\\
To prove Theorem~\ref{q1} below, we will need the following lemma:

\begin{lem}\label{charx}
For $(a_1,a_2,a_3,a_4)\in (K_\Kp^\times)^4$, The following two conditions are equivalent:
\begin{itemize}
    \item [(1)] There exists $(x_1,x_2,x_3,x_4)\in (K_\Kp)^4\backslash\{(0,0,0,0)\}$ such that $a_1x_1^2+a_2x_2^2+a_3x_3^2+a_4x_4^2=0$.
    \item [(2)] There exists nonzero $x'\in K_\Kp$ and $(x_1,x_2,x_3,x_4)\in (K_\Kp)^4\backslash\{(0,0,0,0)\}$ such that 
    $$a_1x_1^2+a_2x_2^2=-a_3x_3^2-a_4x_4^2=x'$$
   
\end{itemize}
\end{lem}
\begin{proof}
(2) $\Rightarrow$ (1) is obvious.\\

For (1) $\Rightarrow$ (2), suppose there exists
\[
(x_1, x_2, x_3, x_4) \in (K_\Kp)^4 \setminus \{(0,0,0,0)\}
\]
such that
\[
a_1 x_1^2 + a_2 x_2^2 + a_3 x_3^2 + a_4 x_4^2 = 0.
\]

If $a_1 x_1^2 + a_2 x_2^2 \neq 0$, take $x' = a_1 x_1^2 + a_2 x_2^2$, and we are done.  
If $a_1 x_1^2 + a_2 x_2^2 = 0$, then the quadratic form $a_1 x_1^2 + a_2 x_2^2$ is isotropic, and hence universal by 42:10 of \cite{iqf}.  
In this case, we can take $x' = -a_3 x_3^2 - a_4 x_4^2$ for any $(x_3, x_4) \in (K_\Kp)^2$ such that $x' \neq 0$.
\end{proof}
\begin{thm} \label{q1}
Let $K$ be a global field with $\textup{char}(K) \neq 2$, and let $f$ be a non-degenerate quadratic form in $m$ variables over $K$. Then the following statements hold:
\begin{itemize}
    \item [(1)] $f$ is isotropic over $K$ if and only if $f$ is isotropic over $K_\Kp$ for every prime $\Kp \in S_K$.
    \item [(2)]If $m \geq 5$ and $\Kp \in S_K$ is a non-archimedean place, then $f$ is isotropic over $K_\Kp$.
    \item [(3)] If $m = 4$, $\Kp \in S_K$ is non-archimedean, and 
    \[
    f = a_1 x_1^2 + a_2 x_2^2 + a_3 x_3^2 + a_4 x_4^2,
    \] 
    then $f$ is anisotropic over $K_\Kp$ if and only if 
    \[
    a_1 a_2 a_3 a_4 \in K_\Kp^{\times 2} \quad \text{and} \quad (a_1, a_2)_\Kp = -(-a_3, -a_4)_\Kp.
    \]
\end{itemize}
\end{thm}
\begin{proof}
The first item follows from Theorem 66:1 in \cite{iqf}, and the second item from Theorem 63:19 in \cite{iqf}.  
The third item is stated in Theorem 6 of Chapter IV in \cite{Serre1973ACI} for $K = \mathbb{Q}$, and the same proof naturally generalizes to arbitrary global fields.  
We provide a detailed proof of this generalization below, following the method in \cite{Serre1973ACI}.\\

By Lemma~\ref{charx}, the form $f = a_1 x_1^2 + a_2 x_2^2 + a_3 x_3^2 + a_4 x_4^2$ is isotropic over $K_\Kp$ if and only if there exists a nonzero $x' \in K_\Kp$ such that the two forms $a_1 x_1^2 + a_2 x_2^2$ and $-a_3 x_3^2 - a_4 x_4^2$ both represent $x'$.  
Let $x$ denote the coset of $x'$ in $K_\Kp^\times / K_\Kp^{\times 2}$.

The statement that $a_1 x_1^2 + a_2 x_2^2$ represents $x'$ is equivalent to 
\[
(x a_1, x a_2)_\Kp = 1.
\]
After a calculation using the properties of the Hilbert symbol, we get
\begin{align*}
(x a_1, x a_2)_\Kp 
&= (x a_1, x)_\Kp , (x a_1, a_2)_\Kp \\
&= (x,x)_\Kp , (a_1,x)_\Kp , (x,a_2)_\Kp , (a_1,a_2)_\Kp \\
&= (x,-1)_\Kp , (x,a_1)_\Kp , (x,a_2)_\Kp , (a_1,a_2)_\Kp \\
&= (x,-a_1 a_2)_\Kp , (a_1,a_2)_\Kp.
\end{align*}

Hence, $a_1 x_1^2 + a_2 x_2^2$ represents $x'$ if and only if
\[
(x,-a_1 a_2)_\Kp = (a_1,a_2)_\Kp,
\]
which is equivalent to
\[
x \in T_{-a_1 a_2}^{-1}((a_1,a_2)_\Kp).
\]

Similarly, $-a_3 x_3^2 - a_4 x_4^2$ represents $x'$ if and only if
\[
(x,-a_3 a_4)_\Kp = (-a_3,-a_4)_\Kp,
\]
which is equivalent to
\[
x \in T_{-a_3 a_4}^{-1}((-a_3,-a_4)_\Kp).
\]

Therefore, the absence of such an $x$ is equivalent to
\[
T_{-a_1 a_2}^{-1}((a_1,a_2)_\Kp) \cap T_{-a_3 a_4}^{-1}((-a_3,-a_4)_\Kp) = \emptyset.
\]
By Lemma~\ref{half}, this is in turn equivalent to
\[
a_1 a_2 = a_3 a_4 \text{ up to a square in } K_\Kp^\times \quad \text{and} \quad (a_1,a_2)_\Kp = -(-a_3,-a_4)_\Kp,
\]
as claimed.
\end{proof}

\subsection{Global class field theory}

In this section, we recall some standard results from global class field theory.  
We begin by introducing the notion of a modulus, which generalizes the role of an integer $m$ in classical modular arithmetic.

\begin{defin}
    Let $K$ be a global field and $\CO_K$ its ring of integers. 
    A \emph{modulus} of $K$ is a formal product of primes of $K$
    \[
        \Km = \prod_{\Kp} \Kp^{\Km(\Kp)},
    \]
    where each exponent $\Km(\Kp) \geq 0$, and the following conditions are satisfied:
    \begin{itemize}
        \item $\Km(\Kp)=0$ for all but finitely many primes $\Kp$;
        \item $\Km(\Kp)=0$ for every infinite complex prime $\Kp$;
        \item $\Km(\Kp)\in\{0,1\}$ for every infinite real prime $\Kp$.
    \end{itemize}
\end{defin}
In the number field case, a modulus can be written as
\[
    \Km = \Km_0 \cdot \Km_\infty,
\] 
where $\Km_0$ denotes the finite part and $\Km_\infty$ the infinite part. 
The notion of congruence extends naturally to this setting.

\begin{defin}\label{decong}
    Let $a,b \in K^\times$, and let $\Km$ be a modulus of $K$. 
    For a finite or real infinite prime $\Kp$, we define
    \[
        a \equiv^* b \pmod{\Kp^{\Km(\Kp)}}
    \]
    as follows:
    \begin{itemize}
        \item If $\Kp$ is finite, then 
        \[
            a\equiv^* b \pmod{\Kp^{\Km(\Kp)}} 
            \quad \Longleftrightarrow \quad 
            v_\Kp\!\left(\tfrac{a}{b}-1\right)\geq \Km(\Kp).
        \]
        \item If $\Kp$ is real, then 
        \[
            a\equiv^* b \pmod{\Kp} 
            \quad \Longleftrightarrow \quad 
            \frac{a}{b}>0
        \]
        under the real embedding associated with $\Kp$.
    \end{itemize}
    Moreover, we write 
    \[
        a \equiv^* b \pmod{\Km}
    \]
    if the above holds for every $\Kp$ with $\Km(\Kp)>0$.
\end{defin}
Let $I_K$ denote the group of fractional ideals of $K$, and let 
\[
    I_\Km \coloneqq \langle \Kp \subseteq \CO_K : \Km(\Kp)=0 \rangle
\]
be the subgroup of $I_K$ generated by all prime ideals $\Kp$ of $K$ that do not divide $\Km$.
\begin{defin}
    The \emph{ray modulo $\Km$} is
    \[
        K_{\Km,1} = \{ a \in K^\times : a \equiv^* 1 \pmod{\Km} \}.
    \]
\end{defin}

Every $\alpha \in K^\times$ defines a principal fractional ideal
\[
    (\alpha) \coloneqq \alpha \CO_K \;\in\; I_K.
\]
If $\alpha \in K_{\Km,1}$, then $(\alpha) \in I_\Km$, so we obtain a well-defined homomorphism
\[
    i:K_{\Km,1}\;\longrightarrow\; I_\Km,\qquad 
    \alpha \longmapsto (\alpha).
\]

\begin{defin}
    Let 
    \[
        P_\Km \coloneqq i(K_{\Km,1}) \;\leq\; I_\Km.
    \]
    The \emph{ray class group modulo $\Km$} is defined as
    \[
        C_\Km \coloneqq I_\Km / P_\Km.
    \]
\end{defin}
Let $L/K$ be a finite abelian extension. 
For a prime ideal $\mathfrak{p}$ of $K$ unramified in $L$, 
denote by
\[
    (\mathfrak{p}, L/K) \;\in\; \mathrm{Gal}(L/K)
\]
the associated Frobenius automorphism.  

If $\Km$ is a modulus of $K$ divisible by all primes that ramify in $L$, 
the \emph{global Artin map} is the homomorphism
\[
    \psi_{L/K} : I_\Km \;\longrightarrow\; \mathrm{Gal}(L/K),
\]
defined on a fractional ideal 
\[
    \mathfrak{a} \;=\; \prod \mathfrak{p}_i^{e_i} \;\in\; I_\Km
\]
by
\[
    \psi_{L/K}(\mathfrak{a}) \;\coloneqq\; \prod (\mathfrak{p}_i, L/K)^{e_i}.
\]
\begin{thm}[Artin Reciprocity]\label{artin}
Let $L/K$ be a finite abelian extension. 
There exists a modulus $\mathfrak{m}$ of $K$, containing all primes of $K$ that ramify in $L$, 
such that the global Artin map
\[
    \psi_{L/K} : I_\mathfrak{m} \longrightarrow \mathrm{Gal}(L/K)
\]
is surjective, with kernel
\[
    \ker(\psi_{L/K}) \;=\; P_\mathfrak{m} \cdot N_{L/K}\big(I_L(\mathfrak{m}')\big),
\]
where $\mathfrak{m}'$ is the modulus of $L$ associated to $\mathfrak{m}$, defined as follows.

If 
\[
    \mathfrak{m} \;=\; 
    \prod_{\mathfrak{p}} \mathfrak{p}^{n_\mathfrak{p}}
    \cdot \prod_{\substack{\omega \text{ real}}} \omega^{n_\omega},
    \qquad n_\omega \in \{0,1\},
\]
then
\[
    \mathfrak{m}' \;=\; 
    \prod_{\mathfrak{P} \mid \mathfrak{p}} \mathfrak{P}^{n_\mathfrak{p}} 
    \cdot \prod_{\substack{\Omega \mid \omega \text{ real}}} \Omega^{n_\omega},
\]
where the first product runs over all finite primes $\mathfrak{P}$ of $L$ above $\mathfrak{p}$, 
and the second product runs over all real embeddings $\Omega$ of $L$ extending the real embedding $\omega$ of $K$.
\end{thm}
We call a modulus $\mathfrak{m}$ as in Theorem \ref{artin} an \emph{admissible modulus} for the extension $L/K$.

Theorem \ref{artin} is stated for number fields and function fields in various standard references. 
For instance, it appears as Theorem~2.1 in \cite{KTU} in the case of global function fields, 
and as Theorem~3.3 in \cite{JP} in the case of number fields. 
For a standard exposition of the Artin symbol, we refer to \cite{Lang1986}, 
where a proof of Artin reciprocity is given in Theorems~1--3.

\section{Main Method}

In this section, we recall the methods used in \cite{10.2307/24735167,JP,KTU,DIT} to construct diophantine sets. 
These methods will serve as the main ingredients in the proof of our main result.

\subsection{Notations and Some Examples}

Let $K$ be a field, and let $A$ be a finite-dimensional central simple algebra over $K$. We define
\[
    S(A/K) \coloneqq \{\textup{Trd}(x) : x \in A, \ \textup{Nrd}(x) = 1\}, \qquad 
    T(A/K) \coloneqq S(A/K) + S(A/K),
\]
where $\textup{Trd}$ and $\textup{Nrd}$ denote the reduced trace and norm, respectively.

From now on, let $K$ be a global field of $\char(K)\neq 2$, and let $A$ be the quaternion algebra $H_{a,b}$ for $a,b\in K^\times$ defined by
\[
    H_{a,b} \coloneqq K \oplus \alpha K \oplus \beta K \oplus \alpha\beta K,
\]
where $\alpha^2 = a$, $\beta^2 = b$, $\alpha\beta = \beta\alpha$.\\
Note that for the quaternion algebra $H_{a,b}$ over a global field $K$, we have
\[
T(A/K) \coloneqq S(A/K) + S(A/K) = S(A/K) - S(A/K),
\]
since for any $x \in A$ with $\mathrm{Nrd}(x)=1$, we have $-x \in A$ and $\mathrm{Nrd}(-x)=\mathrm{Nrd}(x)=1$, 
and moreover $\mathrm{Trd}(-x) = -\mathrm{Trd}(x)$. 
Hence $x \in S(A/K)$ implies $-x \in S(A/K)$, which ensures that $S(A/K)+S(A/K) = S(A/K)-S(A/K)$.

For $\Kp\in S_K^f$, let $R_\Kp$ be its valuation ring in $K_\Kp$, and define
    \[
        \CO_\Kp \coloneqq (\CO_K)_\Kp = R_\Kp \cap K
    \]
    as the local ring of $\Kp$ in $K$.
For simplicity, we write $S(A/K)$ as $S_{a,b}$ and $T(A/K)$ as $T_{a,b}$, since $A$ is the quaternion algebra as stated.
\begin{lem}\label{Tab}
    If $A/K$ splits at all real places of $K$ (this holds if $K$ is a global function field, or if $K$ is a number field and $a,b \in K^{\times}$ are totally positive), then 
    \[
    T_{a,b} = \bigcap_{\mathfrak{q} \in \Delta_{a,b}} \CO_{\mathfrak{q}},
    \]
    where $\Delta_{a,b}$ is the finite set of places of $K$ at which $A$ does not split.
\end{lem}

\begin{proof}
    This is a special case of Proposition 2.9 in \cite{DIT}.
\end{proof}
\begin{rem} \label{uniT}
$T_{a,b}$ is a diophantine set over $K$ by definition.  
This fact forms the basis for proving that many other sets are diophantine in \cite{KTU, DIT, JP}.  
Moreover, $T_{a,b}$ is uniformly existentially defined in $a$ and $b$, since we have the following formula for $S_{a,b}$:
\[
S_{a,b} = \{ 2 x_1 \in K : \exists x_2, x_3, x_4 \text{ with } x_1^2 - a x_2^2 - b x_3^2 + a b x_4^2 = 1 \},
\]
together with the definition $T_{a,b} = S_{a,b} + S_{a,b}$.
\end{rem}
The following definition appears in Definition 9 of \cite{10.2307/24735167} for $K = \BQ$, 
and is also given after the proof of Lemma 3.8 in \cite{JP} as well as after the proof of Lemma 3.3 in \cite{KTU}.

\begin{defin}
Let $L \coloneqq K(\sqrt{a}, \sqrt{b}) / K$ be an abelian extension with Galois group isomorphic to $(\BZ/2\BZ)^2$, 
and let $\psi_{L/K} : I_{\mathfrak{m}} \rightarrow \textup{Gal}(L/K)$ be the global Artin map with admissible modulus $\mathfrak{m}$.  
We define the following notations:

\begin{itemize}
    \item $\BP(p) \coloneqq \{\Kp \in S_K^f : v_{\Kp}(p) \text{ is odd}\}$.
    \item For $(i,j) \in \textup{Gal}(L/K)$, $i,j \in \{\pm 1\}$, we define
    \[
    \BP^{(i,j)} \coloneqq \{\Kp \in S_K^f : \Kp \in I_{\mathfrak{m}} \text{ and } \psi_{L/K}(\Kp) = (i,j)\}.
    \]
    \item $\BP^{(i,j)}(p) \coloneqq \BP(p) \cap \BP^{(i,j)}$.
\end{itemize}

Here, $(1,1)$ is the identity in the Galois group.  
$(-1,1)$ is the element which sends $\sqrt{a}$ to $-\sqrt{a}$ and fixes $\sqrt{b}$.  
$(1,-1)$ is the element which sends $\sqrt{b}$ to $-\sqrt{b}$ and fixes $\sqrt{a}$.  
Finally, $(-1,-1)$ is the element which sends $\sqrt{a}$ to $-\sqrt{a}$ and $\sqrt{b}$ to $-\sqrt{b}$.
\end{defin}

The following lemma is a summary of some results in \cite{JP} and \cite{KTU}.\\
\begin{lem}\label{para}
    (1) If $K$ is a number field, there exist $a, b \in K^{\times}$ with disjoint supports and $\Gal(K(\sqrt{a}, \sqrt{b})/K) \cong (\BZ/2\BZ)^2$ such that for any $p \in K^{\times}$ coprime to the admissible modulus $\mathfrak{m}$, we have:
    \begin{itemize}
        \item $\BP^{(-1,-1)}(p) = \Delta_{a,p} \cap \Delta_{b,p}$,
        \item $\BP^{(-1,1)}(p) = \Delta_{a,p} \cap \Delta_{ab,p}$,
        \item $\BP^{(1,-1)}(p) = \Delta_{b,p} \cap \Delta_{ab,p}$.
    \end{itemize}
    (2) If $K$ is a global function field with $\operatorname{char}(K) \neq 2$, there exist $a, b, c, d \in K^{\times}$ such that $a$ and $b$ have disjoint supports and $\Gal(K(\sqrt{a}, \sqrt{b})/K) \cong (\BZ/2\BZ)^2$. For any $p \in K^{\times}$ coprime to the admissible modulus $\mathfrak{m}$, we have:
    \begin{itemize}
        \item $\BP^{(-1,-1)}(p) = \Delta_{a,p} \cap \Delta_{b,p}$,
        \item $\BP^{(-1,1)}(p) = \Delta_{a,p} \cap \Delta_{ab,p} \cap \Delta_{a,cp}$,
        \item $\BP^{(1,-1)}(p) = \Delta_{b,p} \cap \Delta_{ab,p} \cap \Delta_{b,dp}$.
    \end{itemize}
\end{lem}
\begin{proof}
    We first address (1), which was proved in \cite{JP}. By Lemma 3.9 of \cite{JP}, for any $p$ coprime to $\mathfrak{m}$, the sets on the left and right hand sides of each equation differ at most by the primes dividing the modulus. Then, Lemma 3.19 and Lemma 3.20 in \cite{JP} provide a careful choice of $a$ and $b$ such that these sets are equal.
    
    Part (2) follows from Corollary 3.11 in \cite{KTU}. 
    
    The fact that $a$ and $b$ have disjoint supports is guaranteed by Lemma 3.19 (4) of \cite{JP} and Lemma 3.8 (4) of \cite{KTU}.
\end{proof}
From now on, unless otherwise stated, the parameters $a, b$ (or $a, b, c, d$) are as fixed in Lemma \ref{para}, and $L$ denotes the field $K(\sqrt{a}, \sqrt{b})$, depending on whether $K$ is a number field or a function field.
\begin{defin}
    Let $K$ be a global field with $\operatorname{char}(K) \neq 2$, and keep the constants $a, b, c, d$ as above. For each $\sigma \in \Gal(L/K)$, we define the ring $R_p^{\sigma}$ as follows:
    \begin{enumerate}
        \item[(a)] If $K$ is a number field, then:
        \begin{itemize}
            \item $R_{p}^{(-1,-1)} \coloneqq \bigcap_{\mathfrak{p} \in \Delta_{a,p} \cap \Delta_{b,p}} \CO_{\mathfrak{p}}$
            \item $R_p^{(1,-1)} \coloneqq \bigcap_{\mathfrak{p} \in \Delta_{ab,p} \cap \Delta_{b,p}} \CO_{\mathfrak{p}}$
            \item $R_p^{(-1,1)} \coloneqq \bigcap_{\mathfrak{p} \in \Delta_{a,p} \cap \Delta_{ab,p}} \CO_{\mathfrak{p}}$
            \item $R_{p,q}^{(1,1)} \coloneqq \bigcap_{\mathfrak{p} \in \Delta_{ap,q} \cap \Delta_{bp,q}} \CO_{\mathfrak{p}}$
        \end{itemize}
        \item[(b)] If $K$ is a global function field, then:
        \begin{itemize}
            \item $R_{p}^{(-1,-1)} \coloneqq \bigcap_{\mathfrak{p} \in \Delta_{a,p} \cap \Delta_{b,p}} \CO_{\mathfrak{p}}$
            \item $R_p^{(1,-1)} \coloneqq \bigcap_{\mathfrak{p} \in \Delta_{ab,p} \cap \Delta_{b,p} \cap \Delta_{a,cp}} \CO_{\mathfrak{p}}$
            \item $R_p^{(-1,1)} \coloneqq \bigcap_{\mathfrak{p} \in \Delta_{a,p} \cap \Delta_{ab,p} \cap \Delta_{b,dp}} \CO_{\mathfrak{p}}$
            \item $R_{p,q}^{(1,1)} \coloneqq \bigcap_{\mathfrak{p} \in \Delta_{ap,q} \cap \Delta_{bp,q}} \CO_{\mathfrak{p}}$
        \end{itemize}
    \end{enumerate}
\end{defin}
For $K=\BQ$, the definitions of $R_{p}^\sigma$ and $R_{p,q}^{(1,1)}$ above coincide with Definition 7 in \cite{10.2307/24735167}. Their generalizations to number fields and global fields are given in Definition 3.10 of \cite{JP} and Definition 3.12 of \cite{KTU}, respectively.

The following lemma is well known (see Lemma 3.2.6 and Theorem 3.2.7 in \cite{engler2005valued}).
\begin{lem}\label{ja}
Let $\Delta$ be a finite subset of $S_K^f$. Recall that $J(R)$ denotes the Jacobson radical of a ring $R$.  
If
\[
R=\bigcap_{\Kp \in \Delta} \CO_{\Kp},
\]
then
\[
J(R)=\bigcap_{\Kp \in \Delta} \Kp \CO_{\Kp}.
\]
\end{lem}

The following lemma describes the behavior of addition for finite intersections of valuation rings and their Jacobson radicals.

\begin{lem}\label{sj}
Let $\Delta_1, \Delta_2$ be two finite subsets of $S_K^f$. Then:
\begin{itemize}
    \item[(1)]
    \[
    \bigcap_{\Kp \in \Delta_1} \CO_{\Kp}
    +
    \bigcap_{\Kp \in \Delta_2} \CO_{\Kp}
    =
    \bigcap_{\Kp \in \Delta_1 \cap \Delta_2} \CO_{\Kp}.
    \]
    \item[(2)]
    \[
    \bigcap_{\Kp \in \Delta_1} \Kp \CO_{\Kp}
    +
    \bigcap_{\Kp \in \Delta_2} \Kp \CO_{\Kp}
    =
    \bigcap_{\Kp \in \Delta_1 \cap \Delta_2} \Kp \CO_{\Kp}.
    \]
\end{itemize}
\end{lem}

\begin{proof}
For both (1) and (2), the left-hand sides are obviously contained in the right-hand sides. 
For the reverse inclusion, let $t \in \{0,1\}$ and $x \in K$ be such that $v_\Kp(x) \ge t$ for all $\Kp \in \Delta_1 \cap \Delta_2$.

By weak approximation, there exists $x_1 \in K$ such that
\[
v_\Kp(x_1) \ge t \quad \text{for all } \Kp \in \Delta_1,
\]
and
\[
v_\Kp(x_1 - x) \ge t \quad \text{for all } \Kp \in \Delta_2 \setminus \Delta_1.
\]
Set $x_2 = x - x_1$. Then, for all $\Kp \in \Delta_1 \cap \Delta_2$, we have
\[
v_\Kp(x_2) = v_\Kp(x - x_1) 
\ge \min\{v_\Kp(x), v_\Kp(x_1)\} 
\ge t.
\]
Thus $x = x_1 + x_2$, which is of the form appearing on the left-hand side.
\end{proof}
The following definition summarizes Definition~3.21 and Definition~3.24 in \cite{JP}, which also appear as Definition~3.13 and Definition~3.16 in \cite{KTU}.
\begin{defin}
For each $\sigma\in \Gal(L/K)$, define
\[
\Phi_\sigma
\coloneqq\{p\in K^{\times}:(p)\in I_{\mathfrak{m}},\ \psi_{L/K}((p))=\sigma,\textup{ and }\BP(p)\subset \BP^{(1,1)}\cup \BP^{\sigma}\},
\]
\[
\widetilde{\Phi_\sigma}\coloneqq{K^{\times}}^2\cdot \Phi_\sigma,
\]
and
\[
\Psi_K
\coloneqq\Big\{(p,q)\in \widetilde{\Phi_{(1,1)}}\times \widetilde{\Phi_{(-1,-1)}}\ \Big|\ 
\prod_{\Kp\mid\mathfrak{m}}(ap,q)_{\mathfrak{p}}=-1
\textup{ and }
p\in a\cdot {K^{\times}}^2\cdot (1+J(R_q^{(-1,-1)}))\Big\}.
\]
\end{defin}
For later use, we summarize several results from \cite{KTU} and \cite{JP} as follows.

\begin{pty}\label{kp}
Let $K$ be a global field with $\operatorname{char}(K) \neq 2$, and keep the notations introduced above. Then the following statements hold:
\begin{itemize}
    \item[(1)] The sets
    $\Phi_{\sigma}$, $\Psi_K$, $R_p^{\sigma}$, $J(R_p^{\sigma})$, $R_{p,q}^{(1,1)}$, and $J(R_{p,q}^{(1,1)})$
    are all diophantine sets over $K$.
    
    \item[(2)] For any $p\in \Phi_{\sigma}$ and any $\sigma\in \Gal(L/K)$ with $\sigma\neq (1,1)$, the set $\BP^{\sigma}(p)$ is nonempty.
    
    \item[(3)] Let $\sigma\in \Gal(L/K)$ with $\sigma\neq (1,1)$, and let $\Kp\nmid \Km$ be a prime of $K$ satisfying $\psi_{L/K}(\Kp)=\sigma$. Then there exists an element $p\in \Phi_{\sigma}$ such that
    \[
    \BP^{\sigma}(p)=\{\Kp\}.
    \]
    
    \item[(4)] For any $(p,q)\in \Psi_K$, we have
    \[
    \emptyset\neq \Delta_{ap,q}\cap \Delta_{bp,q}\subset I_{\Km}.
    \]
    
    \item[(5)] For each prime ideal $\Kp_0$ satisfying $\Kp_0\nmid \Km$ and $\psi_{L/K}(\Kp_0)=(1,1)$, there exists $(p,q)\in \Psi_K$ such that
    \[
    \Delta_{ap,q}\cap \Delta_{bp,q}=\{\Kp_0\}.
    \]
    Moreover, we may require that $v_{\Kp_0}(p)=1$, $v_{\Kp_0}(q)=0$, and that $q$ is not a square modulo $\Kp_0$.
\end{itemize}
\end{pty}
\begin{proof}
This proposition follows from Lemma 3.14 and Lemma 3.17 in \cite{KTU} for global function fields, and from Lemma 3.22 and Lemma 3.25 in \cite{JP} for number fields. 
Only the last assertion in (5) requires additional justification, which we provide below.

\medskip
\noindent
\textbf{Number field case.}
The proof of Lemma 3.25(c) in \cite{JP} is based on an explicit construction of the pair $(p,q)$. 
In the choice of $p,q\in K$, condition (B) for $q$ requires that the Legendre symbol $\big(\frac{q}{\Kp_0}\big)$ is equal to $-1$. 
This implies that $v_{\Kp_0}(q)=0$ and that $q$ is not a square modulo $\Kp_0$. 

Condition (4) for $p$ is $(q,p)_{\Kp_0}=-1$. 
Using the explicit formula for the Hilbert symbol, we compute
\begin{align*}
-1=(q,p)_{\Kp_0}
&=(-1)^{v_{\Kp_0}(q)v_{\Kp_0}(p)\frac{|\BF_{\Kp_0}|-1}{2}}
\big(\tfrac{q}{\Kp_0}\big)^{v_{\Kp_0}(p)}
\big(\tfrac{p}{\Kp_0}\big)^{v_{\Kp_0}(q)} \\
&=\big(\tfrac{q}{\Kp_0}\big)^{v_{\Kp_0}(p)},
\end{align*}
since $v_{\Kp_0}(q)=0$. 
As $\big(\tfrac{q}{\Kp_0}\big)=-1$, this implies that $v_{\Kp_0}(p)$ is odd. 
By multiplying $p$ by a square in $K^\times$, we may further assume that $v_{\Kp_0}(p)=1$.

\medskip
\noindent
\textbf{Global function field case.}
This is addressed in the proof of Lemma 3.17(3) in \cite{KTU}. 
In the construction, the element $q$ again satisfies $\big(\tfrac{q}{\Kp_0}\big)=-1$. 
Since $\Kp_0\in \Delta_{ap,q}$, we have $(ap,q)_{\Kp_0}=-1$. 
By computing the Hilbert symbol and using the fact that $\Kp_0\nmid a$, we conclude that $v_{\Kp_0}(p)$ is odd. 
This is stated explicitly in the third-to-last line of the proof of Theorem 1.3 in \cite{KTU}.
\end{proof}
\begin{rem}\label{uniR}
In fact, the sets
$R_p^{\sigma}$, $J(R_p^{\sigma})$, $R_{p,q}^{(1,1)}$, and $J(R_{p,q}^{(1,1)})$
are not only diophantine, but are defined uniformly in the parameters $p$ and $q$.
More precisely, the following sets are diophantine:
\[
\{(p,x)\in K^{\times}\times K : x\in R_p^{\sigma}\},
\]
\[
\{(p,x)\in K^{\times}\times K : x\in J(R_p^{\sigma})\},
\]
\[
\{(p,q,x)\in K^{\times}\times K^{\times}\times K : x\in R_{p,q}^{(1,1)}\},
\]
\[
\{(p,q,x)\in K^{\times}\times K^{\times}\times K : x\in J(R_{p,q}^{(1,1)})\}.
\]

For the rings $R_p^{\sigma}$ and $R_{p,q}^{(1,1)}$, this uniformity follows from the uniform definability of $T_{a,b}$ (see Remark~\ref{uniT}), together with Lemma~\ref{Tab} and Lemma~\ref{sj}(1).
For the corresponding Jacobson radicals, the uniformity follows from Lemma~3.17 of \cite{JP}, Lemma~\ref{ja}, and Lemma~\ref{sj}(2).
\end{rem}

\subsection{Diagonalization of quadratic forms}
Since we always work over a global field $K$ with characteristic not equal to $2$, we obtain the following proposition concerning the diagonalization of quadratic forms.

\begin{pty}\label{diag}
Let
\[
f=\sum_{1\leq i,j\leq m} a_{ij} x_i x_j
\]
be a quadratic form in $m$ variables over $K$. Define the symmetric matrix
$A=[A_{ij}]_{1\leq i,j\leq m}$ by setting $A_{ii}=a_{ii}$ and
$A_{ij}=\frac{a_{ij}+a_{ji}}{2}$ for $i\neq j$.
Then there exist $m^2$ polynomials
\[
(h_{kl})_{1\leq k,l\leq m}\subset \BZ[\overline{t_{ij}},\overline{u_{rs}}],
\]
depending only on $m$, where $\overline{t_{ij}}$ denotes the $m^2$ variables
$(t_{ij})_{1\leq i,j\leq m}$ and $\overline{u_{rs}}$ denotes the $m^2$ variables
$(u_{rs})_{1\leq r,s\leq m}$, such that the following statements are equivalent:

\begin{itemize}
    \item[(1)] The quadratic form $f$ is anisotropic over $K$.
    
    \item[(2)] There exist $b_1,\dots,b_m\in K^{\times}$ such that the diagonal quadratic form
    \[
    \sum_{i=1}^m b_i x_i^2
    \]
    is anisotropic over $K$, and there exist $(u_{rs})_{1\leq r,s\leq m}\in (K)^{m^2}$ satisfying
    \[
    \begin{cases}
        h_{kl}(\overline{A_{ij}},\overline{u_{rs}})=0, & 1\leq k\neq l\leq m, \\
        h_{kk}(\overline{A_{ij}},\overline{u_{rs}})=b_k, & 1\leq k\leq m.
    \end{cases}
    \]
\end{itemize}
\end{pty}
\begin{proof}
Consider the ring of $m\times m$ matrices with entries in $\mathbb{Z}[\overline{t_{ij}}, \overline{u_{rs}}]$, denoted by
\[
M_{m\times m}(\mathbb{Z}[\overline{t_{ij}}, \overline{u_{rs}}]).
\]
Define two matrices
\[
T \coloneqq [t_{ij}]_{1\leq i,j\leq m} \in M_{m\times m}(\mathbb{Z}[\overline{t_{ij}}]) 
\subset M_{m\times m}(\mathbb{Z}[\overline{t_{ij}}, \overline{u_{rs}}]),
\]
\[
U \coloneqq [u_{rs}]_{1\leq r,s\leq m} \in M_{m\times m}(\mathbb{Z}[\overline{u_{rs}}]) 
\subset M_{m\times m}(\mathbb{Z}[\overline{t_{ij}}, \overline{u_{rs}}]).
\]
Then
\[
U^t T U \in M_{m\times m}(\mathbb{Z}[\overline{t_{ij}}, \overline{u_{rs}}]),
\]
and the $(k,l)$-th entry of $U^t T U$ is a polynomial
\[
h_{kl} \in \mathbb{Z}[\overline{t_{ij}}, \overline{u_{rs}}].
\]
By construction, the polynomials $h_{kl}$ depend only on the number of variables $m$.

\medskip
For matrices
\[
A \coloneqq [A_{ij}]_{1\leq i,j\leq m} \in M_{m\times m}(K), \qquad
C \coloneqq [C_{rs}]_{1\leq r,s\leq m} \in M_{m\times m}(K),
\]
we define the evaluation of $H=(h_{kl})$ at $A$ and $C$ by
\[
H(A,C)_{kl} \coloneqq h_{kl}(\overline{A_{ij}}, \overline{C_{rs}}) \in K.
\]
This is well-defined since if $\char(K)=0$, then $\mathbb{Z}\subset K$, and if $\char(K)=p>0$, we interpret the polynomials modulo $p$.

\medskip
Let $B = \text{diag}(b_1,\dots,b_m)$ be a diagonal matrix. Then the matrix equation
\[
C^t A C = B
\]
is equivalent to
\[
\begin{cases}
h_{kl}(\overline{A_{ij}}, \overline{C_{rs}}) = 0, & 1\leq k\neq l \leq m,\\
h_{kk}(\overline{A_{ij}}, \overline{C_{rs}}) = b_k, & 1\leq k \leq m.
\end{cases}
\]

\medskip
If $f$ has full rank $m$, the existence of $b_1,\dots,b_m \in K^\times$ and $C\in (K)^{m^2}$ satisfying this system is equivalent to the diagonalization of $f$ as
\[
\sum_{i=1}^m b_i x_i^2.
\]
Hence $f$ is isotropic over $K$ if and only if the diagonal form $\sum_{i=1}^m b_i x_i^2$ is isotropic over $K$. 

If the rank of $f$ is less than $m$, then $f$ always isotropic, and the equivalence holds trivially.

\medskip
In conclusion, condition (2) implies (1) by the above discussion, and condition (1) implies (2) since the symmetric matrix $A$ is diagonalizable over $K$.
\end{proof}

Under the assumption that our main theorem (Theorem~\ref{t1}) holds, Corollary~\ref{t2} follows immediately from the diagonalization process in Proposition~\ref{diag}. We now give a precise derivation.
\begin{proof}
Let $f = \sum_{1\leq i,j \leq m} a_{ij} x_i x_j$ be a quadratic form defined over $K$. Let the matrix $A=[A_{ij}]_{1\leq i,j\leq m}\in M_{m\times m}(K)$ be the matrix defined in Proposition \ref{diag}.
By Proposition~\ref{diag}, $f$ is isotrpoic over $K$ if and only if there exist
$b_1,\dots,b_m \in K^\times$ and $C_{rs}\in K$ ($1\leq r,s\leq m$) such that
\begin{align*}
h_{kl}(\overline{A_{ij}}, \overline{C_{rs}}) &= 0, & 1\leq k \neq l \leq m,\\
h_{kk}(\overline{A_{ij}}, \overline{C_{rs}}) &= b_k, & 1\leq k \leq m,
\end{align*}
and the diagonal quadratic form $\sum_{i=1}^m b_i x_i^2$
is isotropic over $K$.

By Theorem~\ref{t1}, the condition that the diagonal form $\sum_{i=1}^m b_i x_i^2$ is isotropic over $K$ is diophantine. 
Therefore, the set of coefficients $(a_{ij})_{1\leq i,j\leq m}$ for which $f$ is anisotropic is diophantine, proving Corollary~\ref{t2}.
\end{proof}

\section{Proof of the main results}\subsection{Dealing with a Single Prime}
Our strategy consists of partitioning the set of primes into finitely many classes and handling the primes in each class in a uniform manner. As a preliminary step, we consider the local case for a single prime.

The following lemma, which is a slightly different version of Lemma 4.4 in \cite{KTU}, is restated here for completeness, as it constitutes a key ingredient in our proof.
\begin{lem}
     \label{emlemma}
         Let $\Kp$ be a finite prime of $K$ such that the residue field $\mathbb{F}_{\mathfrak{p}}$ has odd characteristic. Fix $p,s\in K^\times$ such that $v_\Kp(p)$ is odd, $v_\Kp(s)=0$ and the reduction $\textup{red}_\Kp(s)$ is not a square in $\mathbb{F}_{\mathfrak{p}}$. Then for $x,y\in K^\times$, $(x,y)_\Kp=-1$ if and only if 
    $$((x\in p\cdot {K^\times}^2\cdot \CO_\Kp^\times)\wedge(y\textup{ or }-xy\in s\cdot {K^\times}^2\cdot (1+\Kp\CO_\Kp))) $$
    $$\vee ((y\in p\cdot {K^\times}^2\cdot \CO_\Kp^\times)\wedge(x\textup{ or }-xy\in s\cdot {K^\times}^2\cdot (1+\Kp\CO_\Kp))).$$
\end{lem}
\begin{proof}
    We modify the condition $v_{\Kp}(p)=1$ used in Lemma 4.4 of \cite{KTU} to the more general condition that $v_{\Kp}(p)$ is odd. This modification does not alter the set $p \cdot {K^\times}^2 \cdot \CO_{\Kp}^\times$. 
\end{proof}
Besides, we also need the following elementary observation to obtain an equivalent form of Lemma \ref{emlemma}, which is Lemma \ref{eql}.
\begin{lem}\label{bfkey}
    Let $\Kp$ be a finite prime of $K$ such that $|\mathbb{F}_{\Kp}|$ is odd. Fix $p, s \in K^\times$ such that $v_{\Kp}(p) = 1$ and $v_{\Kp}(s) = 0$, with $\operatorname{red}_{\Kp}(s)$ being a non-square in the residue field $\mathbb{F}_{\Kp}$. Then for $x \in K^\times$, the following hold:
    \begin{enumerate}
        \item[(1)] $v_{\Kp}(x)$ is even if and only if $x \in {K^\times}^2 \cdot \CO_{\Kp}^\times$.
        \item[(2)] $v_{\Kp}(x)$ is odd if and only if $x \in p \cdot {K^\times}^2 \cdot \CO_{\Kp}^\times$.
        \item[(3)] $x \in {K_{\Kp}^{\times}}^2$ if and only if $x \in {K^\times}^2 \cdot (1 + \Kp \CO_{\Kp})$.
        \item[(4)] $x \notin {K_{\Kp}^{\times}}^2$ if and only if $x \in s \cdot {K^\times}^2 \cdot (1 + \Kp \CO_{\Kp}) \vee x \in p \cdot {K^\times}^2 \cdot \CO_{\Kp}^\times$.
    \end{enumerate}
    Consequently, we have the following two partitions of $K^\times$:
    \begin{align*}
        K^{\times} &= ( {K^\times}^2 \cdot \CO_{\Kp}^\times ) \sqcup ( p \cdot {K^\times}^2 \cdot \CO_{\Kp}^\times ) \\
        &= ( {K^\times}^2 \cdot (1 + \Kp \CO_{\Kp}) ) \sqcup ( s \cdot {K^\times}^2 \cdot (1 + \Kp \CO_{\Kp}) ) \sqcup ( p \cdot {K^\times}^2 \cdot \CO_{\Kp}^\times ).
    \end{align*}
\end{lem}
\begin{proof}
(1) The valuation $v_{\Kp}(x)$ is even if and only if it can be reduced to $0$ by multiplication by a square in $K^\times$, which is equivalent to $x \in {K^{\times}}^2 \cdot \CO_{\Kp}^{\times}$.

(2) First, if $u \in 1 + \Kp \CO_{\Kp}$, then $u \in {K_{\Kp}^\times}^2$ by the local square theorem (see \cite{iqf}, 63:1). Conversely, suppose $u \in {K_{\Kp}^\times}^2$. By multiplying by a suitable square in $K^\times$, we may assume $v_{\Kp}(u) = 0$. Let $t \in K_{\Kp}^\times$ be such that $u = t^2$. Since $K$ is dense in the local field $K_{\Kp}$, there exists $\eta \in K$ such that $v_{\Kp}(\eta - t) \geq 1$. Given that $v_{\Kp}(t) = 0$, we have $v_{\Kp}(\eta) = \min\{v_{\Kp}(\eta - t), v_{\Kp}(t)\} = 0$, which implies $\eta \in \CO_{\Kp}^\times$. It follows that
\[
u - \eta^2 = (t - \eta)(t + \eta) \in \Kp \CO_{\Kp},
\]
and therefore $u \in \eta^2 + \Kp \CO_{\Kp} = \eta^2(1 + \Kp \CO_{\Kp})$. This shows that $u \in {K^{\times}}^2 \cdot (1 + \Kp \CO_{\Kp})$.

(4) For $x \notin {K_{\Kp}^\times}^2$, we consider two cases based on the parity of the valuation:
\begin{itemize}
    \item If $v_{\Kp}(x)$ is odd, then $x \in p \cdot {K^\times}^2 \cdot \CO_{\Kp}^\times$ follows immediately from (2).
    \item If $v_{\Kp}(x)$ is even, we can write $x = p^{2k}u$ with $v_{\Kp}(u) = 0$ and $u \notin (K_{\Kp}^\times)^2$. According to \cite{iqf}, 63:9, the index $[\CO_{\Kp}^\times : (\CO_{\Kp}^\times)^2]$ is $2$. Since $\operatorname{red}_{\Kp}(s)$ is not a square in the residue field, $us$ must be a square in $K_{\Kp}^\times$. By the result in (2), this implies $sx \in {K^\times}^2 \cdot (1 + \Kp \CO_{\Kp})$, or equivalently, $x \in s \cdot {K^\times}^2 \cdot (1 + \Kp \CO_{\Kp})$.
\end{itemize}
\end{proof}

\begin{lem}
    \label{eql}
    With the same notation as in Lemma \ref{emlemma}, for $x,y\in K^\times$, $(x,y)_\Kp=1$ if and only if 
     $$((x\in {K^\times}^2\cdot \CO_\Kp^\times)\vee(y\textup{ and }-xy\in ({K^\times}^2\cdot (1+\Kp\CO_\Kp))\cup (p\cdot{K^\times}^2\cdot \CO_\Kp^\times)))$$
    $$\wedge ((y\in {K^\times}^2\cdot \CO_\Kp^\times)\vee(x\textup{ and }-xy\in ({K^\times}^2\cdot (1+\Kp\CO_\Kp))\cup (p\cdot {K^\times}^2\cdot \CO_\Kp^\times))).$$
\end{lem}

\begin{proof}
By Lemma \ref{emlemma}, $(x,y)_{\Kp}=-1$ if and only if 
 $$((x\in p\cdot {K^\times}^2\cdot \CO_\Kp^\times)\wedge(y\textup{ or }-xy\in s\cdot {K^\times}^2\cdot (1+\Kp\CO_\Kp))) $$
    $$\vee ((y\in p\cdot {K^\times}^2\cdot \CO_\Kp^\times)\wedge(x\textup{ or }-xy\in s\cdot {K^\times}^2\cdot (1+\Kp\CO_\Kp))).$$
To get $(x,y)_{\Kp}=1$, we take the negation on this condition.
By the partitions established in Lemma \ref{bfkey}, the complement of the set $p \cdot {K^\times}^2 \cdot \CO_{\Kp}^\times$ in $K^\times$ is ${K^\times}^2 \cdot \CO_{\Kp}^\times$. Similarly, the complement of $s \cdot {K^\times}^2 \cdot (1 + \Kp \CO_{\Kp})$ in the union of the two non-square classes is given by
\[
x \in ({K^\times}^2 \cdot (1 + \Kp \CO_{\Kp})) \cup (p \cdot {K^\times}^2 \cdot \CO_{\Kp}^\times).\]
Applying De Morgan's laws and substituting these complements yields the desired characterization.
\end{proof}
In the dyadic case (which means $\Kp|2$), one obtains a characterization of local squares similar to the one presented in Lemma \ref{bfkey} (3). However, due to the requirements of the stronger version of Hensel's Lemma for even primes, a higher power of the prime ideal $\Kp$ is necessary. Specifically, we have the following result:

\begin{lem} \label{dyadic}
    Let $K$ be a number field and $\Kp$ a prime dividing $2$. Let $e$ be the absolute ramification index. Then
    \[
    {K^{\times}}^2 \cdot (1 + \Kp^{2e+1}\CO_{\Kp}) = K \cap (K_{\Kp}^\times)^2.
    \]
\end{lem}
\begin{proof}
Suppose that $u \in 1 + \Kp^{2e+1}\CO_{\Kp}$. By 63:1 of \cite{iqf}, $u$ is a local square, so $u \in {K_{\Kp}^\times}^2$. Since $u \in K$ by assumption, it follows that $u \in K \cap {K_{\Kp}^\times}^2$. This establishes the inclusion ${K^{\times}}^2 \cdot (1 + \Kp^{2e+1}\CO_{\Kp}) \subseteq K \cap {K_{\Kp}^\times}^2$.

Conversely, let $x \in K \cap {K_{\Kp}^\times}^2$. We can write $x = \pi^{2k} u$, where $v_{\Kp}(\pi) = 1$ and $u \in \CO_{\Kp}^\times \cap {K_{\Kp}^\times}^2$. Let $s \in \CO_{\Kp}^\times$ be such that $u = s^2$. If $s \in K$, the result is trivial. If $s \notin K$, then by the density of $K$ in $K_{\Kp}$, there exists $\eta \in K$ such that 
\[ 
v_{\Kp}(\eta - s) \geq 2e + 1. 
\]
Since $v_{\Kp}(s) = 0$ and $2e+1 > 0$, the properties of valuations imply $v_{\Kp}(\eta) = v_{\Kp}(\eta - s + s) = \min\{v_{\Kp}(\eta - s), v_{\Kp}(s)\} = 0$, hence $\eta \in \CO_{\Kp}^\times$. 

To show that $u \eta^{-2} \in 1 + \Kp^{2e+1}\CO_{\Kp}$, we consider the difference:
\[ 
u - \eta^2 = (s - \eta)(s + \eta). 
\]
By our choice of $\eta$, we have $v_{\Kp}(s - \eta) \geq 2e + 1$. For the second factor, since $s + \eta = (s - \eta) + 2\eta$, and noting that $v_{\Kp}(s - \eta) \geq 2e + 1 > e = v_{\Kp}(2\eta)$, we have 
\[ 
v_{\Kp}(s + \eta) = v_{\Kp}(2\eta) = e. 
\]
Therefore, $v_{\Kp}(u - \eta^2) = v_{\Kp}(s - \eta) + v_{\Kp}(s + \eta) \geq (2e + 1) + e = 3e + 1$. 
In particular, this implies $v_{\Kp}(u - \eta^2) \geq 2e + 1$. Since $v_{\Kp}(\eta^2) = 0$, we conclude that
\[ 
u \in \eta^2 (1 + \Kp^{2e+1}\CO_{\Kp}), 
\]
so $x \in {K^{\times}}^2 \cdot (1 + \Kp^{2e+1}\CO_{\Kp})$. This proves the reverse inclusion $K \cap {K_{\Kp}^\times}^2 \subseteq {K^{\times}}^2 \cdot (1 + \Kp^{2e+1}\CO_{\Kp})$.
\end{proof}

\begin{thm}
    \label{one prime}
    Let $\Kp$ be a finite or real infinite prime of $K$. The set $$\{(a_1,a_2,a_3,a_4)\in K^\times\times K^\times\times K^\times\times K^\times:f=a_1x_1^2+a_2x_2^2+a_3x_3^2+a_4x_4^2 \textup{ is anisotropic over }K_\Kp\}$$ is diophantine over $K$.
\end{thm}
\begin{proof}
If $\Kp \in S^f(K)$, then by Theorem \ref{q1}(3), $f$ is anisotropic over $K_{\Kp}$ if and only if $a_1 a_2 a_3 a_4 \in {K_{\Kp}^\times}^2$ and $(a_1, a_2)_{\Kp} = -(-a_3, -a_4)_{\Kp}$. To establish that this condition defines a diophantine set, it suffices to verify that the following three types of sets are diophantine:
\begin{enumerate}
    \item[(1)] $\{x \in K^\times : x \in {K_{\Kp}^\times}^2\}$
    \item[(2)] $\{(x, y) \in K^\times \times K^\times : (x, y)_{\Kp} = 1\}$
    \item[(3)] $\{(x, y) \in K^\times \times K^\times : (x, y)_{\Kp} = -1\}$
\end{enumerate}

Based on the first three chapters of \cite{rumely1980undecidability}, the valuation ring $\CO_{\Kp}$ is diophantine over $K$. Consequently, the maximal ideal $\Kp \CO_{\Kp} = \pi \CO_{\Kp}$ (where $v_{\Kp}(\pi) = 1$) is also diophantine. More precisely, Theorem 1 of \cite{rumely1980undecidability} implies that $\CO_{\Kp}$ is definable, and since the predicate $S_l(x, \vec{c})$ is given by an existential positive formula, it is diophantine.

If $|\mathbb{F}_{\Kp}|$ is odd, the diophantine nature of (2) and (3) follows from Lemma \ref{emlemma} and Lemma \ref{eql}, as these sets are constructed from $\CO_{\Kp}$ and $\Kp \CO_{\Kp}$ using field operations and logical disjunction/conjunction. Additionally, Lemma \ref{bfkey} (3) ensures that the set in (1) is diophantine.

If $K$ is a number field and $\Kp|2$, the index $[K_{\Kp}^\times : {K_{\Kp}^\times}^2]$ is finite (see \cite{iqf}, 63:9). Let $s_1, \ldots, s_n \in K^\times$ be a set of representatives for the cosets of ${K_{\Kp}^\times}^2$ in $K_{\Kp}^\times$. Such representatives can be chosen in $K^\times$ because $K$ is dense in $K_{\Kp}$ and the cosets are open in the local topology. Let $e$ be the absolute ramification index. We define the sets
\[ S_i := s_i \cdot {K^\times}^2 \cdot (1 + \Kp^{2e+1}\CO_{\Kp}) \quad \text{for } i = 1, \ldots, n. \]
These sets $S_1, \ldots, S_n$ form a partition of $K^\times \cap K_{\Kp}^\times$. This follows from two facts:
\begin{itemize}
    \item The sets are disjoint: if $i \neq j$, then $s_i s_j^{-1} \notin {K_{\Kp}^\times}^2$. By Lemma \ref{dyadic}, the set ${K^\times}^2 \cdot (1 + \Kp^{2e+1}\CO_{\Kp})$, being a subset of ${K_{\Kp}^\times}^2$, cannot contain $s_i s_j^{-1}$.
    \item They cover the relevant elements: suppose $x \in K \cap s_i {K_{\Kp}^\times}^2$, then $s_i^{-1}x \in K \cap {K_{\Kp}^\times}^2$. By Lemma \ref{dyadic}, we have $s_i^{-1}x \in {K^\times}^2 \cdot (1 + \Kp^{2e+1}\CO_{\Kp})$, which implies $x \in S_i$.
\end{itemize}
Under this partition, $x \in {K_{\Kp}^\times}^2$ if and only if $x \in S_i$ for $s_i \in {K_{\Kp}^\times}^2$. Furthermore, the Hilbert symbol conditions satisfy:
\[ (x, y)_{\Kp} = 1 \iff (x, y) \in \bigcup_{(s_i, s_j)_{\Kp} = 1} (S_i \times S_j), \]
\[ (x, y)_{\Kp} = -1 \iff (x, y) \in \bigcup_{(s_i, s_j)_{\Kp} = -1} (S_i \times S_j). \]
Since $\Kp^{2e+1}\CO_{\Kp}$ is diophantine, each $S_i$ is diophantine, and so are their finite unions.

For a real infinite prime $\Kp$, let $\omega \colon K \hookrightarrow \mathbb{R}$ be the corresponding embedding. Lemma 4.1 of \cite{KTU} states that the set $P_\omega := \{x \in K^\times : \omega(x) > 0\}$ is diophantine over $K$. The form $f$ is anisotropic over $K_{\Kp} \cong \mathbb{R}$ if and only if it is definite, i.e.,
\[ \left( \bigwedge_{i=1}^4 a_i \in P_\omega \right) \vee \left( \bigwedge_{i=1}^4 -a_i \in P_\omega \right). \]
Since $P_\omega$ is diophantine, this condition defines a diophantine set.
\end{proof}

\subsection{The proof of main theorem}
Now we prove our main result, Theorem~\ref{t1}.

\begin{lem}\label{keylemma}
    The set $$\{(a_1,a_2,a_3,a_4)\in K^\times \times K^{\times} \times K^{\times} \times K^{\times}: f=a_1x_1^2+a_2x_2^2+a_3x_3^2+a_4x_4^2 \textup{ is anisotropic over some }\Kp\in  S_K\}$$ is a diophantine set. 
\end{lem}
\begin{proof}
    For $\sigma \in \Gal(K(\sqrt{a},\sqrt{b})/K)$ with $\sigma \neq (1,1)$, define
\[
s_{\sigma} :=
\begin{cases}
a, & \text{if } \sigma = (-1,\pm 1),\\
b, & \text{if } \sigma = (1,-1).
\end{cases}
\]

We claim that the quadratic form
\[
f = a_1 x_1^2 + a_2 x_2^2 + a_3 x_3^2 + a_4 x_4^2
\]
is anisotropic over some $\Kp \in S_K$ if and only if \emph{one of the following three conditions holds}.

\begin{itemize}

\item[(i)]
There exists a prime $\Kp \mid \Km$ such that $f$ is anisotropic over $K_{\Kp}$.

\item[(ii)]
There exists $\sigma \neq (1,1)$ and an element $p \in \Phi_{\sigma}$ such that
\[
\Bigl(\dagger^{+}_{a_1,a_2} \wedge \dagger^{-}_{-a_3,-a_4}\Bigr)
\;\vee\;
\Bigl(\dagger^{-}_{a_1,a_2} \wedge \dagger^{+}_{-a_3,-a_4}\Bigr)
\;\wedge\; \mathfrak e,
\]
where
\[
\mathfrak e \;:\; a_1 a_2 a_3 a_4 \in K^{\times 2} \cdot (1 + J(R_p^{\sigma})).
\]

Here the formulas $\dagger^{\pm}_{x,y}$ are defined as follows.

\[
\dagger^{-}_{x,y} :=
\begin{aligned}[t]
&\bigl(x \in p \cdot K^{\times 2} \cdot (R_p^{\sigma})^{\times}
\;\wedge\;
(y \text{ or } -xy \in s_{\sigma} \cdot K^{\times 2} \cdot (1 + J(R_p^{\sigma})))\bigr)\\
&\;\vee\;
\bigl(y \in p \cdot K^{\times 2} \cdot (R_p^{\sigma})^{\times}
\;\wedge\;
(x \text{ or } -yx \in s_{\sigma} \cdot K^{\times 2} \cdot (1 + J(R_p^{\sigma})))\bigr),
\end{aligned}
\]

\[
\dagger^{+}_{x,y} :=
\begin{aligned}[t]
&\bigl(x \in K^{\times 2} \cdot (R_p^{\sigma})^{\times}
\;\vee\;
(y \text{ and } -xy \in p \cdot K^{\times 2} \cdot (R_p^{\sigma})^{\times}
\cup K^{\times 2} \cdot (1 + J(R_p^{\sigma})))\bigr)\\
&\;\wedge\;
\bigl(y \in K^{\times 2} \cdot (R_p^{\sigma})^{\times}
\;\vee\;
(x \text{ and } -yx \in p \cdot K^{\times 2} \cdot (R_p^{\sigma})^{\times}
\cup K^{\times 2} \cdot (1 + J(R_p^{\sigma})))\bigr).
\end{aligned}
\]

\item[(iii)]
There exists $(p,q) \in \Psi_K$ such that $q \in (R_{p,q}^{(1,1)})^{\times}$ and
\[
\Bigl(\ddagger^{+}_{a_1,a_2} \wedge \ddagger^{-}_{-a_3,-a_4}\Bigr)
\;\vee\;
\Bigl(\ddagger^{-}_{a_1,a_2} \wedge \ddagger^{+}_{-a_3,-a_4}\Bigr)
\;\wedge\; \mathfrak e',
\]
where $\mathfrak e'$ is defined analogously to $\mathfrak e$, with $R_p^{\sigma}$ replaced by $R_{p,q}^{(1,1)}$.

The formulas $\ddagger^{\pm}_{x,y}$ are given by

\[
\ddagger^{-}_{x,y} :=
\begin{aligned}[t]
&\bigl(x \in p \cdot K^{\times 2} \cdot (R_{p,q}^{(1,1)})^{\times}
\;\wedge\;
(y \text{ or } -xy \in q \cdot K^{\times 2} \cdot (1 + J(R_{p,q}^{(1,1)})))\bigr)\\
&\;\vee\;
\bigl(y \in p \cdot K^{\times 2} \cdot (R_{p,q}^{(1,1)})^{\times}
\;\wedge\;
(x \text{ or } -yx \in q \cdot K^{\times 2} \cdot (1 + J(R_{p,q}^{(1,1)})))\bigr),
\end{aligned}
\]

\[
\ddagger^{+}_{x,y} :=
\begin{aligned}[t]
&\bigl(x \in K^{\times 2} \cdot (R_{p,q}^{(1,1)})^{\times}
\;\vee\;
(y \text{ and } -xy \in p \cdot K^{\times 2} \cdot (R_{p,q}^{(1,1)})^{\times}
\cup K^{\times 2} \cdot (1 + J(R_{p,q}^{(1,1)})))\bigr)\\
&\;\wedge\;
\bigl(y \in K^{\times 2} \cdot (R_{p,q}^{(1,1)})^{\times}
\;\vee\;
(x \text{ and } -yx \in p \cdot K^{\times 2} \cdot (R_{p,q}^{(1,1)})^{\times}
\cup K^{\times 2} \cdot (1 + J(R_{p,q}^{(1,1)})))\bigr).
\end{aligned}
\]
\end{itemize}
The first condition defines a diophantine set, since the set $\Km$ contains only finitely many primes and by Theorem~\ref{one prime}. 
The second and third conditions also define diophantine sets by Proposition~\ref{kp}(1) together with Remark~\ref{uniR}. 
Consequently, the set
\[
\bigl\{(a_1,a_2,a_3,a_4)\in (K^\times)^4 \mid 
f=a_1x_1^2+a_2x_2^2+a_3x_3^2+a_4x_4^2 
\textup{ is anisotropic over some } \Kp\in S_K\bigr\}
\]
is a diophantine set, being a finite union of the three sets defined above.\\
It remains to verify that the claim holds. 
For brevity, we write 
\[
(\dagger_{a_1,a_2}^{+}\wedge \dagger_{-a_3,-a_4}^{-})\vee (\dagger_{a_1,a_2}^{-} \wedge \dagger_{-a_3,-a_4}^{+}) \quad \text{as } \dagger
\]
and 
\[
(\ddagger_{a_1,a_2}^{+}\wedge \ddagger_{-a_3,-a_4}^{-})\vee (\ddagger_{a_1,a_2}^{-} \wedge \ddagger_{-a_3,-a_4}^{+}) \quad \text{as } \ddagger.
\]
The primes $\Kp\nmid \Km$ can be classified according to 
$\sigma \in \Gal(K(\sqrt{a},\sqrt{b})/K)$. 
By Theorem~\ref{q1}(3), the condition that $f$ is anisotropic over $\Kp\nmid \Km$ is equivalent to 
\[
(a_1,a_2)_{\Kp} \neq -(-a_3,-a_4)_{\Kp} \quad \text{and} \quad a_1a_2a_3a_4 \in (K_{\Kp}^\times)^2.
\]
Hence, it suffices to establish the following two equivalences between the statements in our claim and this local condition for the different types of primes $\Kp$.
\begin{itemize}
    \item For $\sigma \neq (1,1)$, the sentence
    \[
    \exists p \in \Phi_{\sigma} \text{ such that } \dagger \wedge \mathfrak{e}
    \]
    is equivalent to the existence of a prime $\Kp$ with $\psi_{L/K}(\Kp) = \sigma$ such that
    \[
    (a_1,a_2)_{\Kp} \neq -(-a_3,-a_4)_{\Kp} \quad \text{and} \quad a_1a_2a_3a_4 \in (K_{\Kp}^\times)^2.
    \]

    \item The sentence
    \[
    \exists (p,q) \in \Psi_K \text{ such that } q \in (R_{p,q}^{(1,1)})^{\times} \text{ and } \ddagger \wedge \mathfrak{e}'
    \]
    is equivalent to the existence of a prime $\Kp$ with $\psi_{L/K}(\Kp) = (1,1)$ such that
    \[
    (a_1,a_2)_{\Kp} \neq -(-a_3,-a_4)_{\Kp} \quad \text{and} \quad a_1a_2a_3a_4 \in (K_{\Kp}^\times)^2.
    \]
\end{itemize}
For $\sigma \neq (1,1)$, suppose there exists $p \in \Phi_{\sigma}$ such that $\dagger$ holds. Then, by Proposition \ref{kp}(2), the set $\BP^{\sigma}(p)$ is nonempty.  

Let $\Kp \in \BP^{\sigma}(p)$. Recall that the formula $\dagger_{x,y}^-$ is given by
\begin{align*}
((x \in p \cdot K^{\times 2} \cdot (R_p^{\sigma})^{\times}) \wedge (y \text{ or } -xy \in s_{\sigma} \cdot K^{\times 2} \cdot (1 + J(R_p^{\sigma})))) \\
\vee ((y \in p \cdot K^{\times 2} \cdot (R_p^{\sigma})^{\times}) \wedge (x \text{ or } -yx \in s_{\sigma} \cdot K^{\times 2} \cdot (1 + J(R_p^{\sigma})))).
\end{align*}

By Lemma \ref{para} and the definition of $R_p^\sigma$, we have $R_p^\sigma \subset \CO_{\Kp}$. Moreover, by Proposition \ref{ja}, $J(R_p^\sigma) \subset \Kp \CO_{\Kp}$.  

Hence, $\dagger_{x,y}^-$ implies $(x,y)_{\Kp} = -1$, because $v_{\Kp}(p)$ is odd, by Lemma \ref{emlemma}, and since $s_\sigma$ is a $\Kp$-adic unit which is not a square (by the construction of $s_\sigma$ and the definition of the Artin map).  
For example, if $\sigma = (1,-1)$, then $b$ is not a square modulo $\Kp$ since $\psi_{L/K}(\Kp) = (1,-1)$.  

By the same reasoning, $\dagger_{x,y}^{+}$ implies $(x,y)_{\Kp} = 1$.  

Finally, the formula $\mathfrak{e}$ implies
\[
a_1 a_2 a_3 a_4 \in K^{\times 2} \cdot (1 + \Kp \CO_{\Kp}),
\]
which means that $a_1 a_2 a_3 a_4$ is a square in $K_{\Kp}$ by Lemma \ref{bfkey}(3).  

Therefore, $\dagger \wedge \mathfrak{e}$ implies
\[
(a_1, a_2)_{\Kp} \neq -(-a_3, -a_4)_{\Kp} \quad \text{and} \quad a_1 a_2 a_3 a_4 \in (K_{\Kp}^\times)^2.
\]
Conversely, suppose there exists a prime $\Kp \nmid \Km$ with $\psi_{L/K}(\Kp) = \sigma$ such that
\[
(a_1, a_2)_{\Kp} \neq -(-a_3, -a_4)_{\Kp}.
\]

Then, by Proposition \ref{kp}(3), we can find $p \in \Phi_{\sigma}$ such that $\BP^{\sigma}(p) = \{\Kp\}$.  
In this situation, for $\sigma \neq (1,1)$, we have
\[
R_p^{\sigma} = \CO_{\Kp}.
\]

Therefore, the formula $\dagger_{x,y}^{\pm}$ can be deduced directly from $(x,y)_{\Kp} = \pm 1$,  
and the formula $\mathfrak{e}$ can be deduced from 
\[
a_1 a_2 a_3 a_4 \in (K_{\Kp}^\times)^2.
\]

Hence, we conclude that
\[
\dagger \wedge \mathfrak{e}
\]
follows from
\[
(a_1, a_2)_{\Kp} \neq -(-a_3, -a_4)_{\Kp} \quad \text{and} \quad a_1 a_2 a_3 a_4 \in (K_{\Kp}^\times)^2.
\]
Conversely, suppose there exists a prime $\Kp \nmid \Km$ with $\psi_{L/K}(\Kp) = \sigma$ such that
\[
(a_1, a_2)_{\Kp} \neq -(-a_3, -a_4)_{\Kp}.
\]

Then, by Proposition \ref{kp}(3), we can find $p \in \Phi_{\sigma}$ such that $\BP^{\sigma}(p) = \{\Kp\}$.  
In this situation, for $\sigma \neq (1,1)$, we have
\[
R_p^{\sigma} = \CO_{\Kp}.
\]

Therefore, the formula $\dagger_{x,y}^{\pm}$ can be deduced directly from $(x,y)_{\Kp} = \pm 1$,  
and the formula $\mathfrak{e}$ can be deduced from 
\[
a_1 a_2 a_3 a_4 \in (K_{\Kp}^\times)^2.
\]

Hence, we conclude that
\[
\dagger \wedge \mathfrak{e}
\]
follows from
\[
(a_1, a_2)_{\Kp} \neq -(-a_3, -a_4)_{\Kp} \quad \text{and} \quad a_1 a_2 a_3 a_4 \in (K_{\Kp}^\times)^2.
\]
A similar argument establishes the second equivalence, using Proposition \ref{kp}(4) and (5).  
We now describe the proof for the case $\psi_{L/K}(\Kp) = (1,1)$.

Assume $\Kp_0 \in I_\Km$ with $\Psi_{L/K}(\Kp_0) = (1,1)$, then by Proposition \ref{kp}, we can find $(p,q) \in \Psi_K$ such that
\[
\Delta_{ap,q} \cap \Delta_{bp,q} = \{\Kp_0\}.
\]
Moreover, we have
\[
v_{\Kp_0}(p) = 1, \quad v_{\Kp_0}(q) = 0,
\]
and $q$ is not a square modulo $\Kp_0$.  
Then, by Lemma \ref{emlemma}, Lemma \ref{bfkey}, and Lemma \ref{eql}, the condition
\[
(a_1,a_2)_{\Kp_0} = -(-a_3,-a_4)_{\Kp_0} \quad \text{and} \quad a_1 a_2 a_3 a_4 \in (K_{\Kp_0}^\times)^2
\]
implies
\[
\ddagger \wedge \mathfrak{e},
\]
since $R_{p,q}^{(1,1)} = \CO_{\Kp_0}$ and $J(R_{p,q}^{(1,1)}) = \Kp_0 \CO_{\Kp_0}$.

Conversely, suppose $(p,q) \in \Psi_K$, then by Proposition \ref{kp}, we have $\Delta_{ap,q} \cap \Delta_{bp,q} \neq \emptyset$.  
For any $\Kp \in \Delta_{ap,q} \cap \Delta_{bp,q}$, since $q \in (R_{p,q}^{(1,1)})^\times$ and 
\[
(ap,q)_\Kp = (bp,q)_\Kp = -1,
\]
we obtain
\[
\left(\frac{q}{\Kp}\right)^{v_{\Kp}(bp)} = \left(\frac{q}{\Kp}\right)^{v_{\Kp}(ap)} = -1.
\]
Hence $q$ is not a square modulo $\Kp$, and $v_{\Kp}(ap), v_{\Kp}(bp)$ are odd.  
Since the supports of $a$ and $b$ are disjoint (Lemma \ref{para}), it follows that $v_{\Kp}(p)$ is odd.  

Again, by Lemma \ref{emlemma}, Lemma \ref{bfkey}, and Lemma \ref{eql}, we conclude that
\[
\ddagger \wedge \mathfrak{e} \implies (a_1,a_2)_\Kp = -(-a_3,-a_4)_\Kp \quad \text{and} \quad a_1 a_2 a_3 a_4 \in (K_\Kp^\times)^2,
\]
since $R_{p,q}^{(1,1)} \subset \CO_{\Kp}$ and $J(R_{p,q}^{(1,1)}) \subset \Kp \CO_{\Kp}$.
\end{proof}

\begin{proof}[proof of  Theorem \ref{t1}: ]
    We proceed by considering each possible number of variables $m$ separately.
    \begin{itemize}
        \item $m=1$. For all $a_1\neq 0$, we have $f=a_1x_1^2$ is anisotropic, and the set $ K^{\times}$
        is obviously diophantine.
        \item $m=2$. For all $a_1,a_2\neq 0$, $f=a_1x_1^2+a_2x_2^2$ is isotropic if and only if $-a_1a_2\in {K^{\times}}^2$, and the set $$\{(a_1,a_2)\in {K^{\times}}\times K^{\times}: -a_1a_2\notin {K^{\times}}^2\}$$ is a diophantine set over $K$ by Theorem 1.4 of \cite{KTU}.
        \item $m=3$. In this case $f=a_1x_1^2+a_2x_2^2+a_3x_3^3$. Then $-a_1a_2\in {K^\times}^2$ means exactly $f=0$ has a nontrivial solution for $x_3=0$. The condition $f=0$ has a nontrivial solution for $x_3\neq 0$ is equivalent to the solubility of the equation $-a_3=a_1x_1^2+a_2x_2^2$, which means exactly 
        $$-a_3=a_1(x_1^2+a_1a_2(\frac{x_2}{a_1})^2)\in a_1\textup{N}_{K(\sqrt{-a_1a_2})/K}(K(\sqrt{-a_1a_2})^{\times}).$$ 
        Hence for all $a_1,a_2,a_3\neq 0$, $f=a_1x_1^2+a_2x_2^2+a_3x_3^3$ is isotropic if and only if $-a_1a_2\in {K^{\times}}^2$ or $-a_3\in a_1\cdot \textup{N}_{K(\sqrt{-a_1a_2})/K}(K(\sqrt{-a_1a_2})^{\times})$.\\
        
        The sets 
        $$\{(a_1,a_2,a_3)\in {K^{\times}}\times K^\times\times K^\times: -a_1a_2\notin {K^{\times}}^2)\}$$ 
        and 
        $$\{(a_1,a_2,a_3)\in {K^{\times}}\times K^\times\times K^\times: -a_3\notin a_1\cdot \textup{N}_{K(\sqrt{-a_1a_2})/K}(K(\sqrt{-a_1a_2})^{\times})\}$$ are all diophantine by Theorem 1.3 of \cite{KTU} and Theorem 1.4 of \cite{KTU}.
        \item $m=4$. This follows from Lemma \ref{keylemma} and Theorem \ref{q1}(1).
        \item $m\geq 5$. By (1) and (2) of Theorem \ref{q1}, in this situation $f$ is anisotropic over $K$ if and only if $f$ is anisotropic over some real infinite primes since $f$ is always isotropic over complex primes and finite primes (63:19 of \cite{iqf}). For a single real infinite prime $\omega$, $f$ is anisotropic over $\omega$ if and only if $$({\wedge_{i=1}^m\omega(a_i)>0})\vee(\wedge_{i=1}^m\omega(-a_i)>0),$$ which gives a diophantine set by Lemma 4.1 of \cite{KTU}. Since there are only finitely many real primes, the set 
        $$\{(a_1,...,a_m)\in(K^\times)^m:\Sigma_ia_ix_i^2 \textup{ is anisotropic over } K\}$$
        is diophantine.
    \end{itemize}
\end{proof}

\section*{Acknowledgments}
I am grateful to my PhD advisor Yongqi Liang for leading me to the question studied in this paper and for many helpful suggestions. I would like to express my sincere thanks to Arno Fehm, my host supervisor during the Sino-German (CSC-DAAD) Postdoc Scholarship Program, who spent a great deal of time carefully reading the manuscript and provided extensive and insightful advice on both the mathematics and the exposition. I would also like to thank Yang Cao and Yongqi Liang for their valuable help and support throughout the application process for the Sino-German (CSC-DAAD) Postdoc Scholarship Program. This work was carried out during a research visit to TU Dresden supported by the Sino-German (CSC-DAAD) Postdoc Scholarship Program and was also supported by the Shandong Provincial Natural Science Foundation (Project No. ZR2024QA184).

\bibliographystyle{alpha}
\bibliography{name}

\end{document}